\documentclass[11pt]{amsart}

\usepackage[letterpaper,margin=0.92in]{geometry}
\usepackage[T1]{fontenc}
\usepackage{mathpazo}
\usepackage{courier}
\usepackage{microtype}
\usepackage{booktabs,tabularx,array}
\usepackage{enumitem}
\usepackage{fancyhdr}
\usepackage{pdfpages}
\usepackage{hyperref}
\usepackage{bookmark}

\hypersetup{
  unicode=true,
  breaklinks=true,
  colorlinks=false,
  pdfborder={0 0 0},
  pdfdisplaydoctitle=true,
  bookmarksopen=true,
  bookmarksdepth=2,
  pdftitle={Resolution of Singularities in Positive Characteristic: Frobenius--Hasse Towers and Exceptional-History Descent. Parts I--IX},
  pdfauthor={Chenxiao Tian},
  pdfsubject={Strong embedded resolution and principalization over perfect fields of positive characteristic},
  pdfkeywords={resolution of singularities, positive characteristic, Frobenius--Hasse tower, exceptional history, principalization, functorial resolution}
}

\newcommand{\SeriesTitle}{Resolution of Singularities in Positive Characteristic}
\newcommand{\SeriesSubtitle}{Frobenius--Hasse Towers and Exceptional-History Descent}
\newcommand{\FrontHeading}[1]{%
  \begin{center}\large\scshape #1\end{center}\vspace{6pt}}
\newcommand{\PartEntry}[4]{%
  \textbf{#1} & \textbf{#2} & #3 & \textbf{\mbox{#4}}\\[2pt]
}

\begin{document}


\thispagestyle{empty}
\pdfbookmark[0]{Title Page}{front-title}
\vspace*{0.55in}
\begin{center}
{\small\scshape A Treatise in Nine Parts}\par
\vspace{0.34in}
{\fontsize{24}{28}\selectfont\bfseries \SeriesTitle\par}
\vspace{0.20in}
{\fontsize{17}{21}\selectfont \SeriesSubtitle\par}
\vspace{0.32in}
\rule{0.72\textwidth}{0.45pt}\par
\vspace{0.31in}
{\Large\bfseries Parts I--IX\par}
\vspace{0.74in}
{\Large Chenxiao Tian\par}
\vspace{0.12in}
{\normalsize Princeton University\par}
\vspace{0.08in}
{\normalsize Princeton, New Jersey 08544, USA\par}
\vfill
{\itshape The singularity may return; its history cannot.\par
Where numerical descent fails, historical descent begins.\par}
\vfill
{\small \texttt{ct3471@princeton.edu}\quad\textbullet\quad
\texttt{ct3471@alumni.princeton.edu}\par}
\vspace{0.18in}
{\small August 2026\par}
\end{center}
\clearpage

\pdfbookmark[0]{Global Abstract}{global-abstract}
\FrontHeading{Global Abstract}

Let $k$ be a perfect field of characteristic $p>0$.  This work constructs a
canonical strong embedded resolution over $k$ by a finite sequence of ordinary
blowups with regular permissible centres, preserving an ordered
simple-normal-crossings boundary at every stage.  The construction is
presentation independent, stable under re-embedding, and functorial for open,
smooth, and \(\acute{e}\)tale pullback and for extension of the perfect ground
field.  It yields strong principalization of coherent ideals, reduced embedded
resolution, and intrinsic resolution by the standard embedding and descent
procedures.

The local theory replaces a marked Rees algebra by its
differential--integral saturation.  Total Hasse operators, coefficient cubes,
and filtered Rees complexes retain the information lost by numerical order in
positive characteristic.  Their transformation under permissible blowup is
encoded by semilinear Frobenius--Hasse sources and by a finite exceptional
ancestry.  Local defect complexes are strictified and routed to surface,
toroidal--monomial, binomial, or additive-type backends; the resulting words
are then serialized on a canonical global nerve.

Termination is expressed at the level of addressed occurrences rather than
by a pointwise scalar invariant.  A global replacement certificate transports
the owner, parent, quotient, trace, and reopening data across macroblocks.  Its
six realization components are supplied by rigid constructor generation,
cross-generation unique addressing, complete wild-capacity control, a
centre-or-typed-exit alternative, structured cofibres, displayed-parent
allocation, and literal terminal truth.  Each completed block therefore
replaces a nonempty multiset of active parent occurrences by strict descendants
in a well-founded dependent order.  The process terminates, and the terminal
truth comparison identifies the exhausted state with the required geometric
normal form.

\vfill
\begin{minipage}{\textwidth}

\end{minipage}
\clearpage

\pdfbookmark[0]{Architecture and Reading Guide}{reading-guide}
\FrontHeading{Architecture and Reading Guide}

The nine parts form a single proof in four layers.  Parts I--III construct and
transport the local semantic packet.  Parts IV--V close its defect calculus and
serialize the resulting local words globally.  Part VI isolates the exact
replacement mechanism responsible for termination and reconstruction.  Parts
VII--IX realize that mechanism object by object and assemble the final
resolution theorem.

\vspace{4pt}
{\small
\begin{tabularx}{\textwidth}{@{}>{\raggedright\arraybackslash}p{0.08\textwidth}
  >{\raggedright\arraybackslash}p{0.31\textwidth}
  >{\raggedright\arraybackslash}X@{}}
\toprule
\textbf{Part} & \textbf{Mathematical layer} & \textbf{Principal output}\\
\midrule
I & Differential--integral local theory & Saturated marked Rees algebra,
total-Hasse activity, coefficient cubes, and filtered transform interface.\\
II & Semilinear height theory & Primitive semilinear sources, height
filtrations, gauge descent, and the prepared generic-entry packet.\\
III & Exceptional transformation & All-chart comparison, literal derived
source, finite-word heredity, endpoint package, and ancestry filtration.\\
IV & Defect closure and backends & Six-row defect ledger, strictification,
surface and monomial backends, additive torsors, and paid handoffs.\\
V & Global serialization & Clean portfolios, canonical refinements, descent
centres, global scheduler, and owner no-reset.\\
VI & Replacement and termination & Dependent occurrence order, renewal forest,
strict endpoint replacement, multiset termination, and reconstruction.\\
VII & Certificate realization I & Realization site, obstruction calculus,
address and quotient criteria, reopening calculus, and chamber assembly.\\
VIII & Certificate realization II & One-step cleavage, prefix heredity,
fresh-run comparison, cross-generation addressing, and wild capacity.\\
IX & Canonical A--F realization & Saturated generation, simultaneous
ancestral-carrier descent, literal quotients, unique addressing, wild
capacity, universal progress, terminal truth, and final assembly.\\
\bottomrule
\end{tabularx}}

\vfill
\textbf{Citation form.}
Individual results are cited by part, section, and local theorem number; for
example, ``Part VI, Theorem 6.7.''  Printed folios and PDF page labels are
part-local, as are theorem and equation numbers.
\clearpage

\pdfbookmark[0]{Global Contents}{global-contents}
\FrontHeading{Global Contents}

{\small
\begin{tabularx}{\textwidth}{@{}>{\raggedright\arraybackslash}p{0.06\textwidth}
  >{\raggedright\arraybackslash}p{0.34\textwidth}
  >{\raggedright\arraybackslash}X
  >{\raggedleft\arraybackslash}p{0.155\textwidth}@{}}
\toprule
\textbf{Part} & \textbf{Title} & \textbf{Scope} & \textbf{Pages}\\
\midrule
\PartEntry{I}{Differential--Integral Local Theory}
  {Saturation, Hasse packets, coefficient cubes, and filtered transforms.}{I-1--I-60}
\PartEntry{II}{Semilinear Height Filtrations and Generic Entry}
  {Primitive sources, height filtrations, gauge torsors, and preparation.}{II-1--II-80}
\PartEntry{III}{Exceptional Transform and Finite-Word Heredity}
  {All-chart transformation, endpoint comparison, and exceptional ancestry.}{III-1--III-70}
\PartEntry{IV}{Defect Closure and Decorated Backends}
  {Defect strictification, certified backends, additive torsors, and payment.}{IV-1--IV-60}
\PartEntry{V}{Global Serialization, Principalization, and Descent}
  {Canonical global nerve, portfolios, scheduler, and owner transport.}{V-1--V-60}
\PartEntry{VI}{Strict Endpoint Replacement, Termination, and Terminal Reconstruction}
  {Replacement certificates, dependent descent, termination, and reconstruction.}{VI-1--VI-50}
\PartEntry{VII}{Certificate Realization I: Obstruction Calculus and Verified Chambers}
  {Realization fields, obstruction calculus, reopening, and chamber criteria.}{VII-1--VII-139}
\PartEntry{VIII}{Certificate Realization II: One-Step Heredity, Cross-Generation Addressing, and Wild Capacity}
  {Cleavage systems, hereditary prefixes, addressing, and complete wild packets.}{VIII-1--VIII-118}
\PartEntry{IX}{Canonical Realization of the A--F Certificate: Saturated Generation, Simultaneous Carrier Descent, and Universal Progress}
  {Foundation theorems, simultaneous ancestral carriers, Targets D--A--E--B--C--F, canonical certificate assembly, and the resolution theorem.}{IX-1--IX-156}
\bottomrule
\end{tabularx}}

\vfill
\begin{center}
\begin{minipage}{0.86\textwidth}
\small\centering
All nine parts begin on recto pages.  The four-page preliminary matter is
numbered in lower-case roman numerals.  Every mathematical Part retains its
own printed folios.
\end{minipage}
\end{center}
\clearpage


\includepdf[pages=-,fitpaper=true,pagecommand={},addtotoc={1,section,0,{Part I: Differential--Integral Local Theory},part-I}]{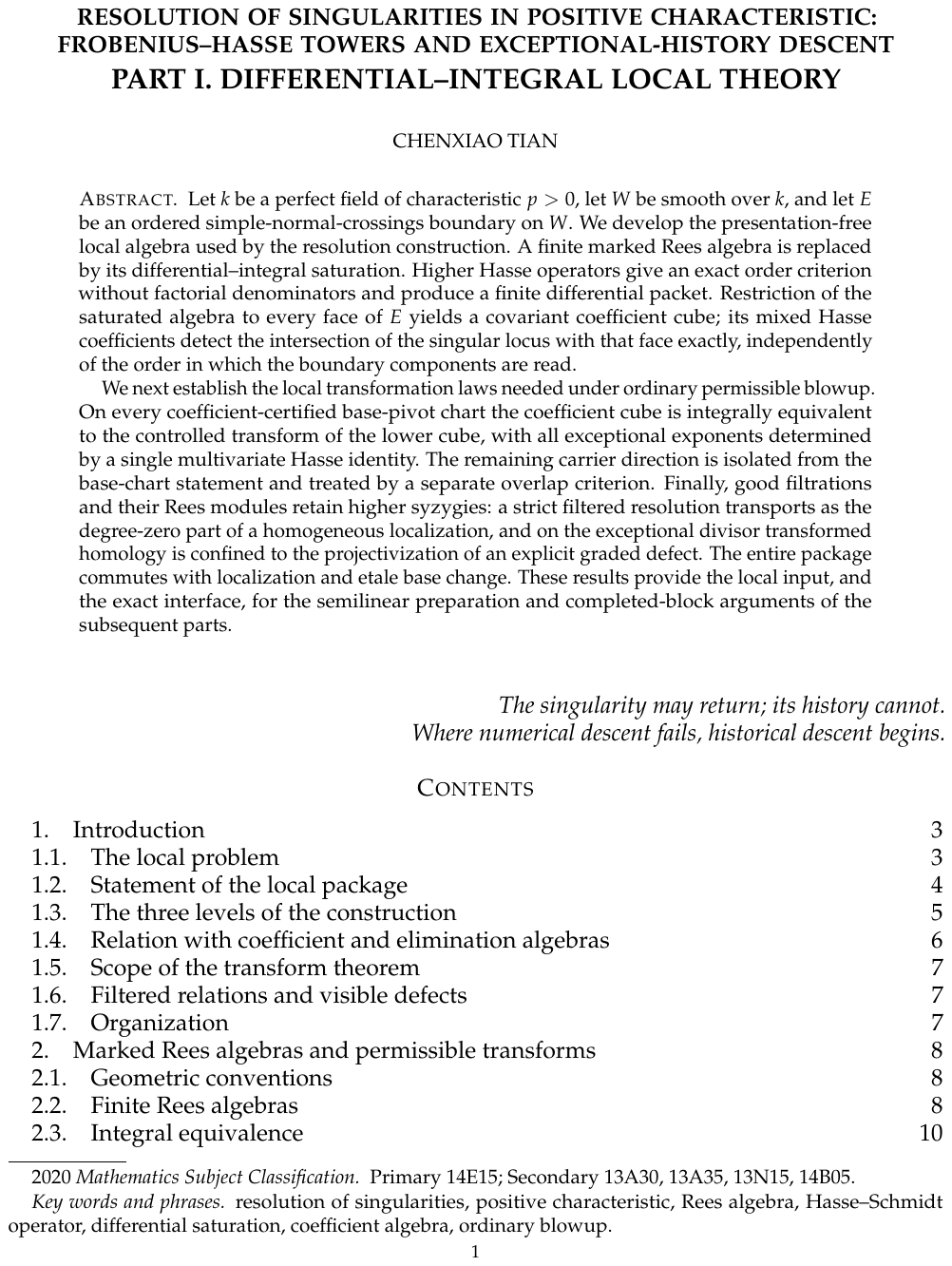}
\includepdf[pages=-,fitpaper=true,pagecommand={},addtotoc={1,section,0,{Part II: Semilinear Height Filtrations and Generic Entry},part-II}]{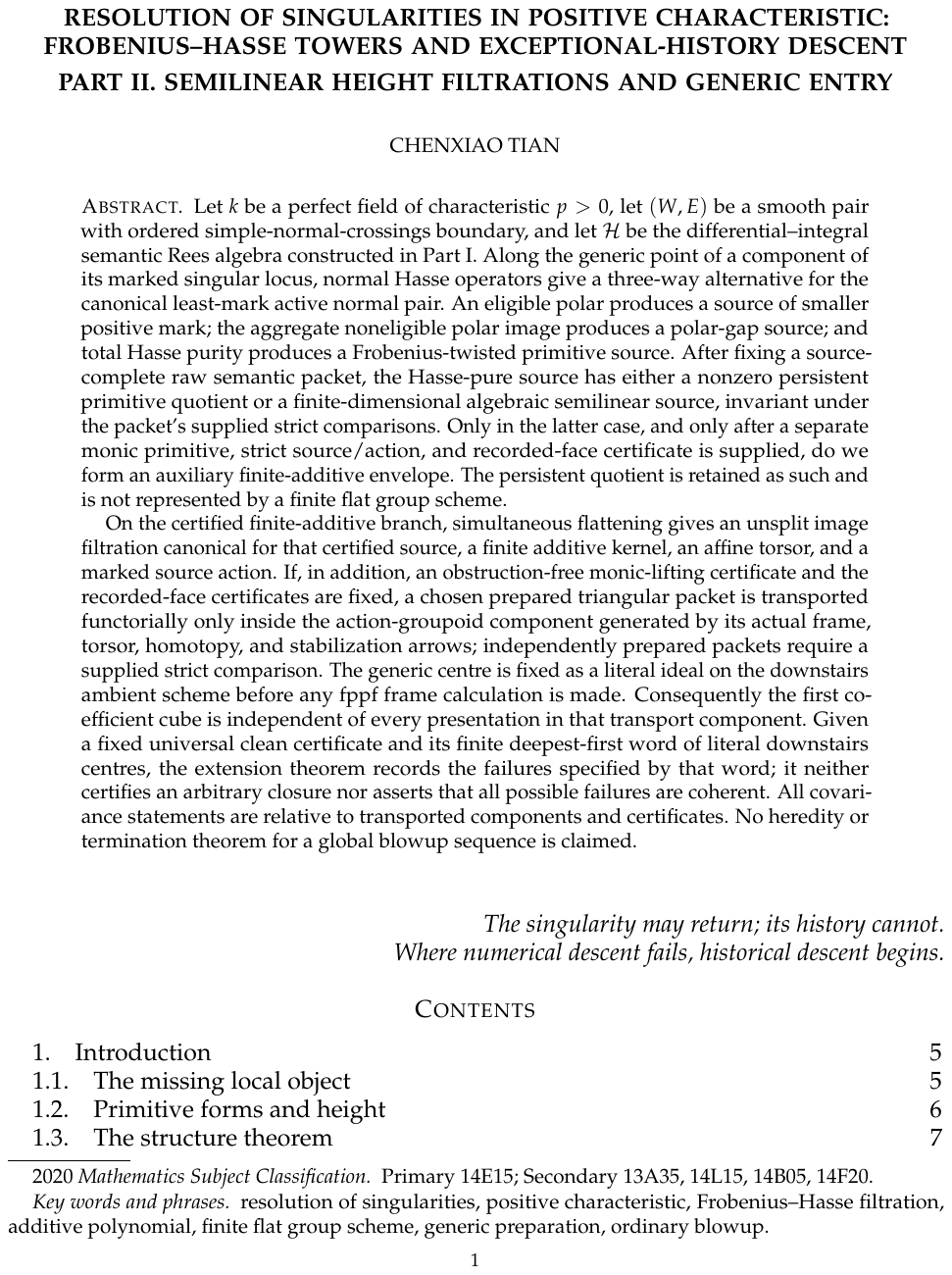}
\includepdf[pages=-,fitpaper=true,pagecommand={},addtotoc={1,section,0,{Part III: Exceptional Transform and Finite-Word Heredity},part-III}]{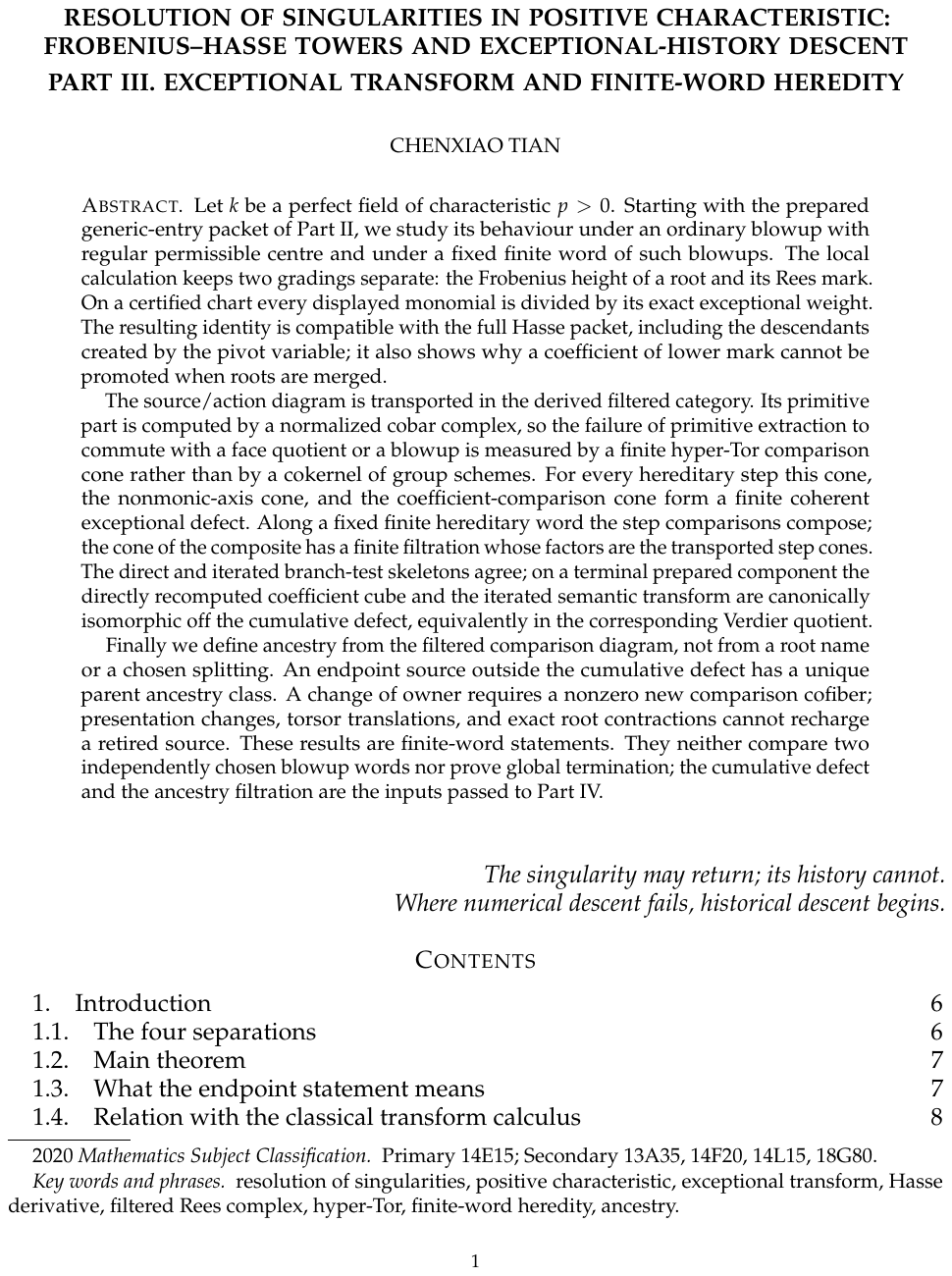}
\includepdf[pages=-,fitpaper=true,pagecommand={},addtotoc={1,section,0,{Part IV: Defect Closure and Decorated Backends},part-IV}]{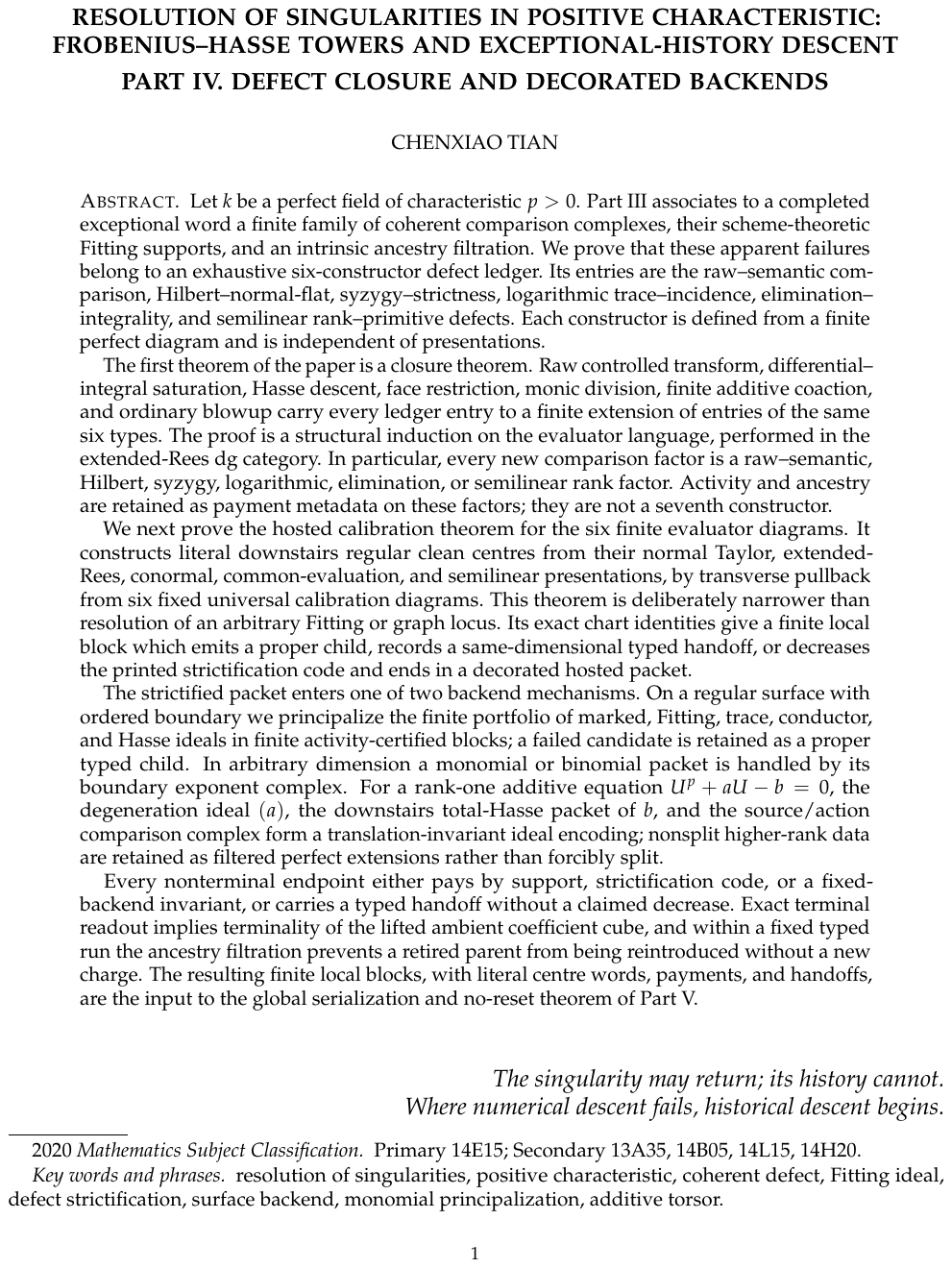}
\includepdf[pages=-,fitpaper=true,pagecommand={},addtotoc={1,section,0,{Part V: Global Serialization, Principalization, and Descent},part-V}]{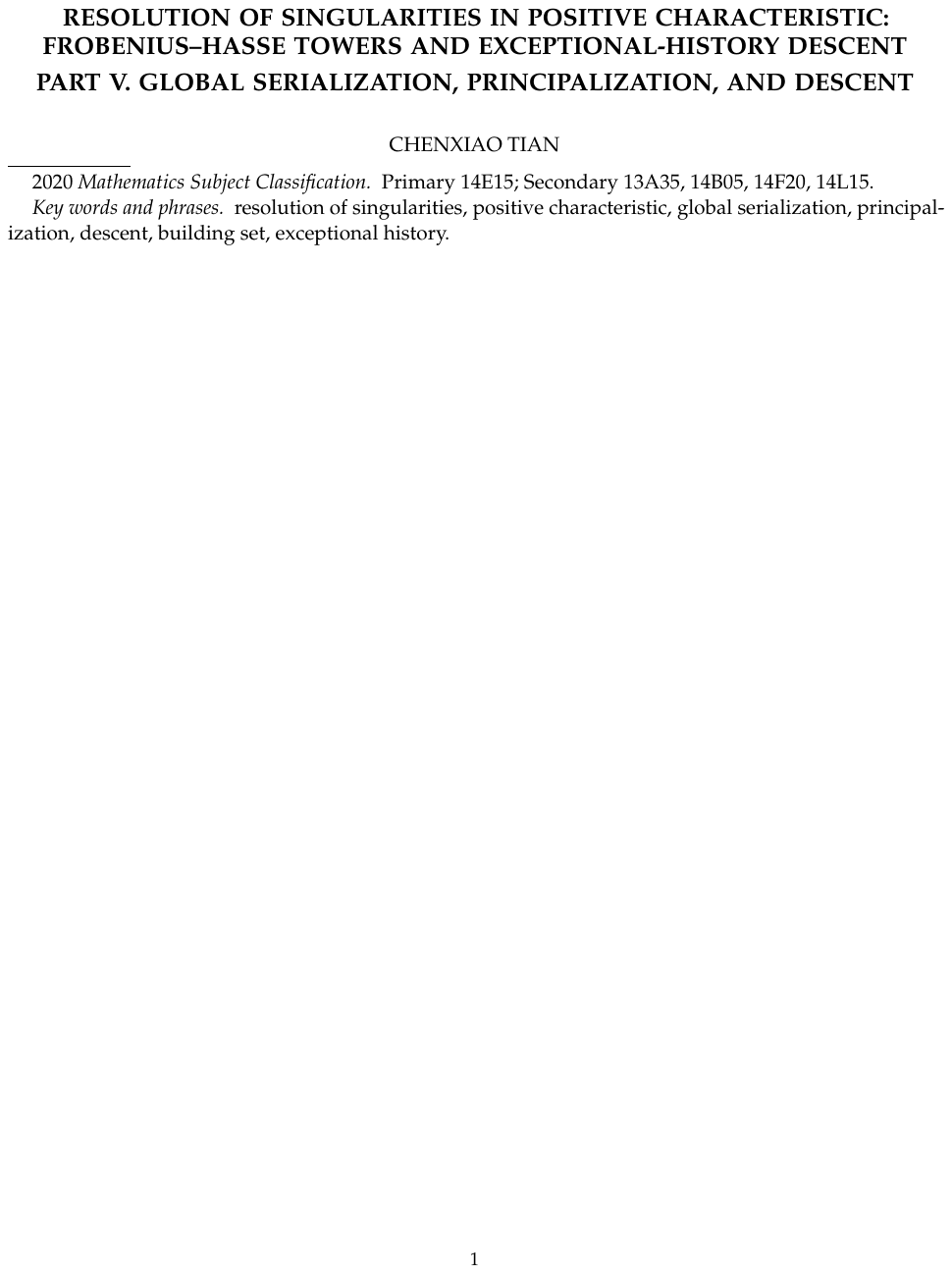}
\includepdf[pages=-,fitpaper=true,pagecommand={},addtotoc={1,section,0,{Part VI: Strict Endpoint Replacement, Termination, and Terminal Reconstruction},part-VI}]{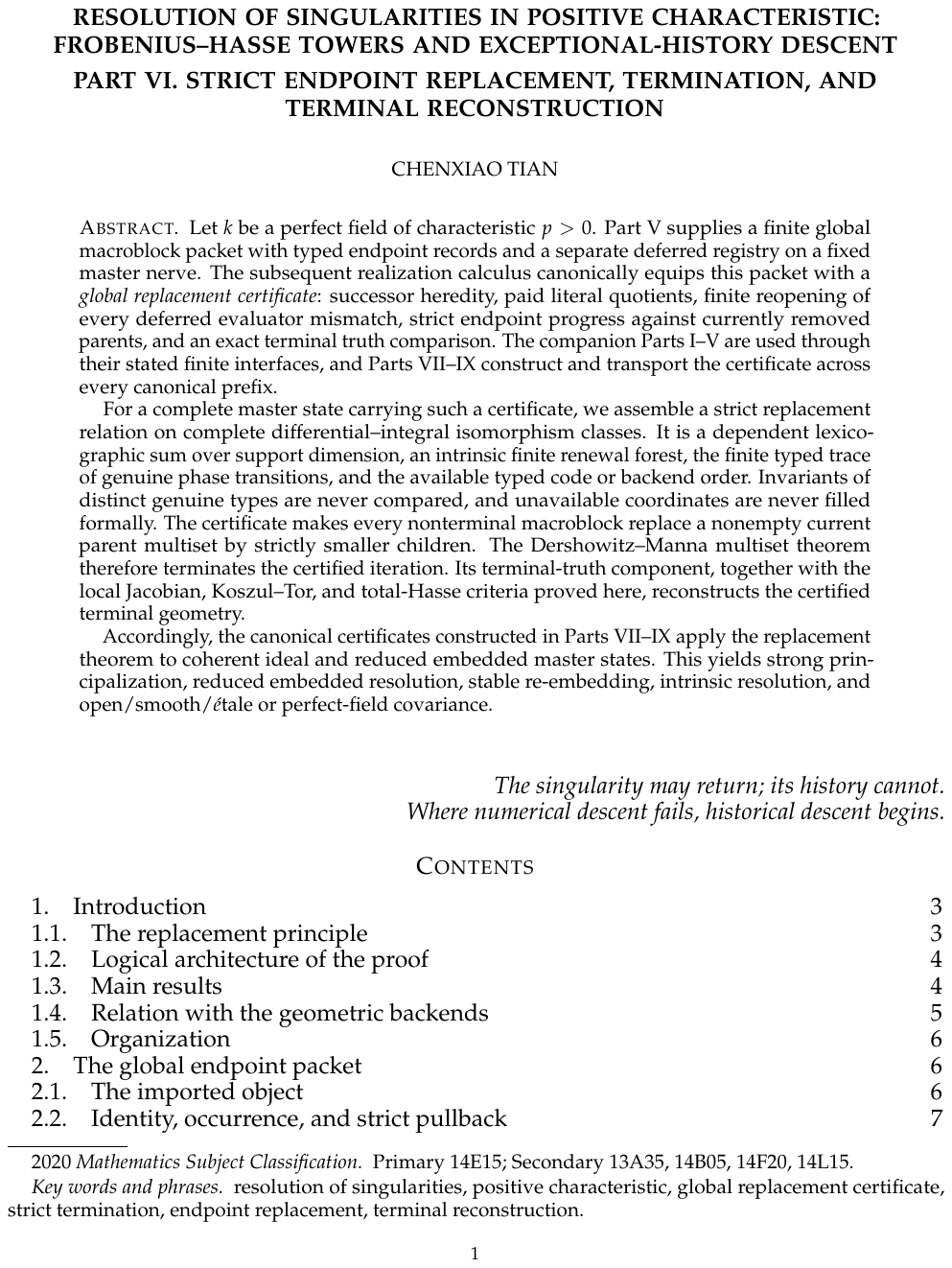}
\includepdf[pages=-,fitpaper=true,pagecommand={},addtotoc={1,section,0,{Part VII: Certificate Realization I},part-VII}]{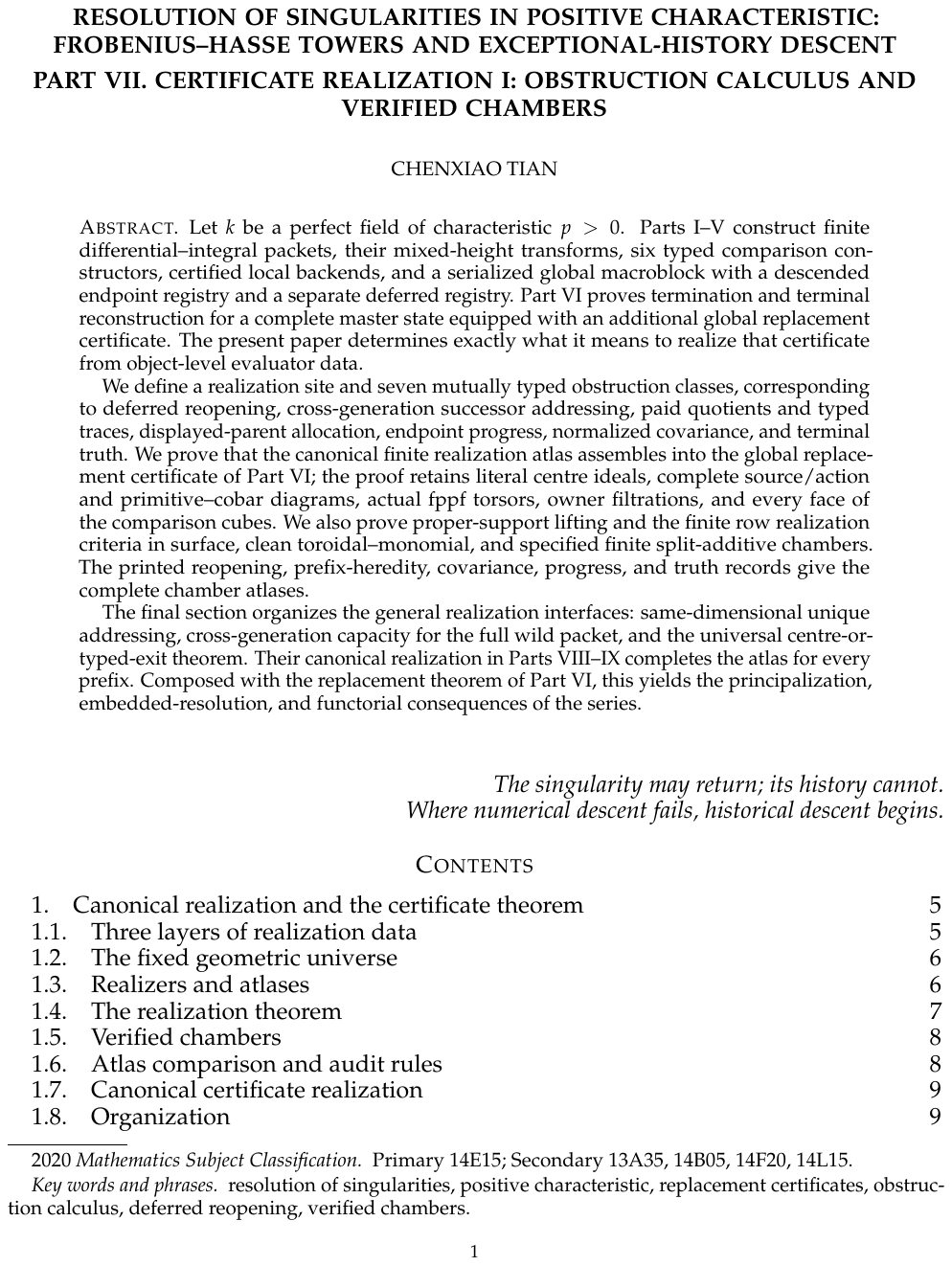}

\clearpage
\thispagestyle{empty}
\null
\clearpage

\includepdf[pages=-,fitpaper=true,pagecommand={},addtotoc={1,section,0,{Part VIII: Certificate Realization II},part-VIII}]{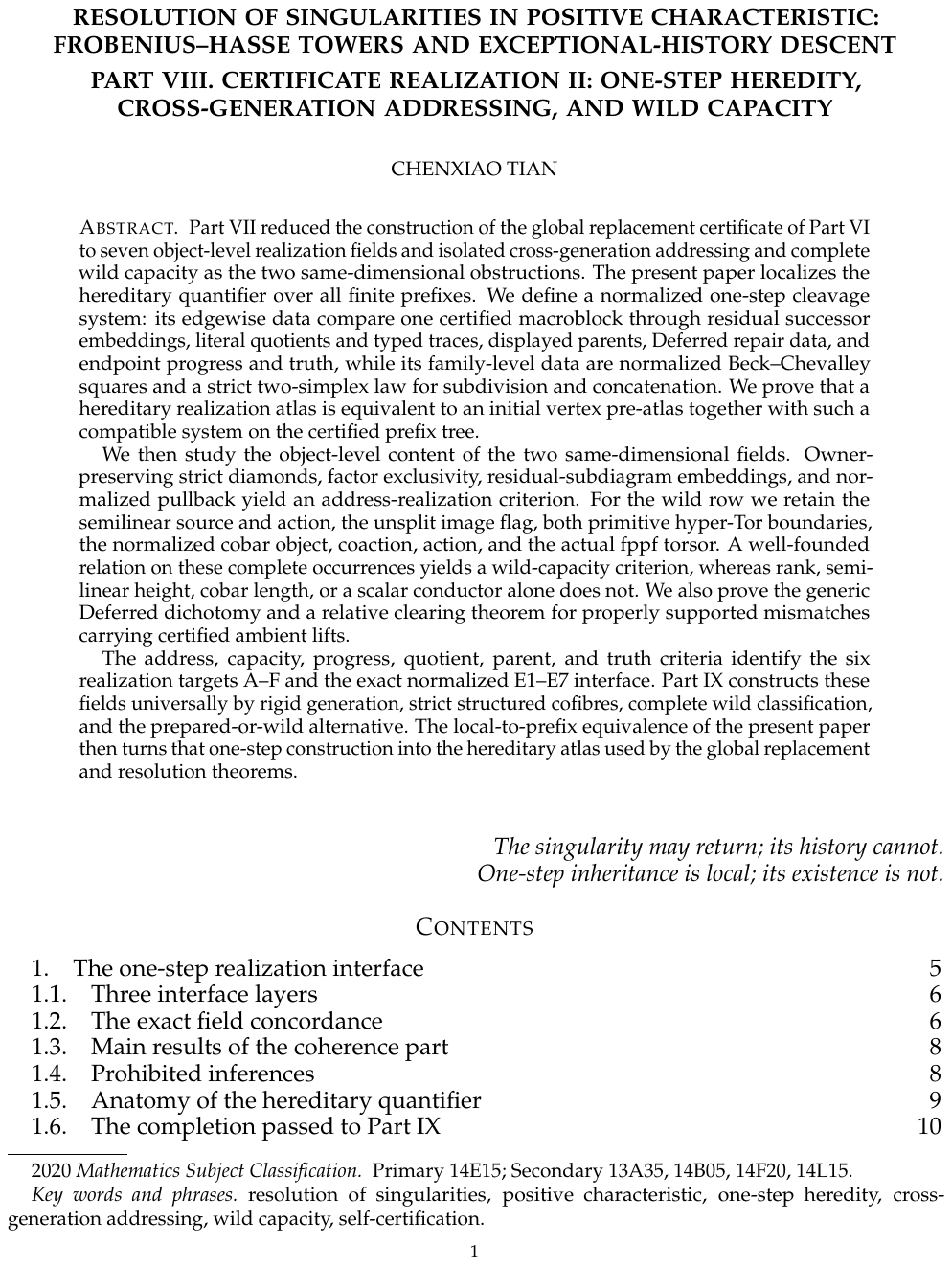}
\includepdf[pages=-,fitpaper=true,pagecommand={},addtotoc={1,section,0,{Part IX: Canonical Realization of the A--F Certificate},part-IX}]{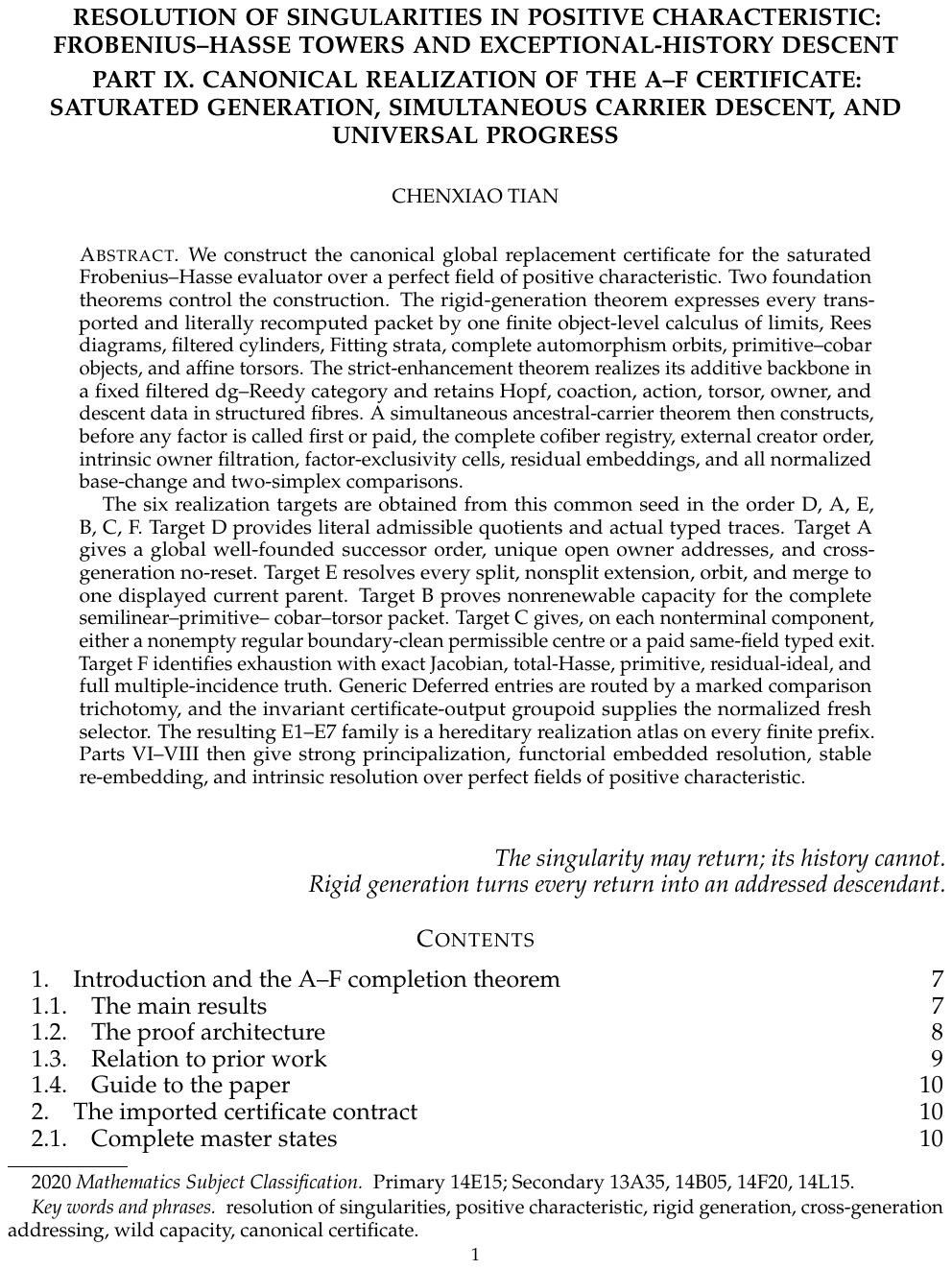}

\end{document}